\UseRawInputEncoding
\documentclass{amsart}
\usepackage{amsfonts}
\newtheorem{theorem}{Theorem}
\theoremstyle{plain}

\newtheorem{corollary}{Corollary}

\newtheorem{definition}{Definition}

\newtheorem{lemma}{Lemma}

\newtheorem{proposition}{Proposition}
\newtheorem{remark}{Remark}

\numberwithin{equation}{section}

\newcommand{\s}{\section}
\newcommand{\R}{\mathbb R}
\newcommand{\C}{\mathbb C}
\newcommand{\lab}{\label}
\newcommand{\bt}{\begin{theorem}}
\newcommand{\et}{\end{theorem}}
\newcommand{\bl}{\begin{lemma}}
\newcommand{\el}{\end{lemma}}
\newcommand{\bd}{\begin{definition}}
\newcommand{\ed}{\end{definition}}
\newcommand{\bc}{\begin{corollary}}
\newcommand{\ec}{\end{corollary}}
\newcommand{\bp}{\begin{proof}}
\newcommand{\ep}{\end{proof}}
\newcommand{\bo}{\begin{proposition}}
\newcommand{\eo}{\end{proposition}}
\newcommand{\br}{\begin{remark}}
\newcommand{\er}{\end{remark}}
\newcommand{\beq}{\begin{equation}}
\newcommand{\eeq}{\end{equation}}
\newcommand{\N}{{\mathbb N}}

\newenvironment{altproof}[1]
{\noindent
{\sl Proof of {#1}}.}
{\nopagebreak\mbox{}\hfill $\Box$\par\addvspace{0.5cm}}

\def\C{\mathbb{C}}
\def\N{\mathbb{N}}

\def\R{\mathbb{R}}

\def\bd{\mathrm{bd}\,}

\begin{document}
\title[Ground state normalized solution]{A new deduce of the strict binding inequality and its application: Ground state normalized solution to Schr\"odinger equations with potential}
\author{X. Zhong}
\address[X. Zhong]{South China Research Center for Applied Mathematics and Interdisciplinary Studies,South China Normal University, Guangzhou 510631, China.}
\address[W. Zou]{Department of Mathematical Sciences, Tsinghua University, Beijing 100084, China.}
\email[X. Zhong]{zhongxuexiu1989@163.com}
\author{W. Zou}
\email[W. Zou]{zou-wm@mail.tsinghua.edu.cn}
\thanks{Supported by NSFC (11801581,11025106, 11371212, 11271386), Guangdong Basic and Applied Basic Research Foundation (2021A1515010034),Guangzhou Basic and Applied Basic Research Foundation(202102020225), Province Natural Science Fund of Guangdong (2018A030310082), the Both-Side NCNU Fund and Tsinghua Fund.}
\date{June 9, 2021}
\subjclass[2010]{35Q55, 35Q51, 35B09, 35C08, 35J20}
\keywords{Schr\"odinger equation, Ground state normalized solution, Global minimizer,Potential.}

\begin{abstract}
In the present paper, we prove the existence of solutions $(\lambda, u)\in \R\times H^1(\R^N)$ to the following elliptic equations with potential
$\displaystyle
-\Delta u+(V(x)+\lambda)u=g(u)\;\hbox{in}\;\R^N,
$
satisfying the normalization constraint $\displaystyle \int_{\R^N}u^2=a>0,$  which is deduced by searching for solitary wave solution to the time-dependent nonlinear Schr\"odinger equations. Besides the importance in the applications, not negligible
reasons of our interest for such problems with potential $V(x)$ are their stimulating and challenging mathematical difficulties. We develop an interesting way basing on iteration and give a new proof of the so-called ``sub-additive inequality", which can simply the standard process in the traditional sense.  Under some very relax assumption on the potential $V(x)$ and some other suitable assumptions on $g$, we can obtain the existence of ground state solution for prescribed $a>0$.
\end{abstract}

\maketitle
\tableofcontents
\s{Introduction}
\renewcommand{\theequation}{1.\arabic{equation}}
\renewcommand{\theremark}{1.\arabic{remark}}
\renewcommand{\thedefinition}{1.\arabic{definition}}

Normalized solutions to semilinear elliptic problems  have been  investigated in many different models, one of the well-established motivation comes from the time-dependent nonlinear Schr\"odinger equations(NLSE). For concreteness, in the present paper we consider the following NLSE
\beq\lab{eq:NLSE}
\begin{cases}
&-i\frac{\partial}{\partial t}\Phi=\Delta \Phi -V(x)\Phi+f(|\Phi|)\Phi=0, \;\;\;\;\;  (x,t)\in \R^N\times \R,\\
&\Phi=\Phi(x,t)\in \C, N\geq 1.
\end{cases}
\eeq
Then the formula  $\Phi(x,t)=e^{i\lambda t}u(x)$   for the  solitary wave solution  establishes the connection between \eqref{eq:NLSE} and  the following kind of elliptic equations:
\beq\lab{eq:NLSES}
-\Delta u+\left(V(x)+\lambda\right) u=g(u)\;\hbox{in}\;\R^N,\;
\eeq
One can choose a fixed $\lambda\in \R$ and search for solutions $u$ of \eqref{eq:NLSES}. We call it the  {\it fixed frequency problem. } There are many mathematical theories and tools  applied to study  it. For example, one can apply the variational method, looking for  the critical points of the associated energy functional
$$J_\lambda[u]:=\frac{1}{2}\int_{\R^N}|\nabla u|^2dx+\frac{1}{2}\int_{\R^N}\left(V(x)+\lambda\right)u^2dx-\int_{\R^N}G(u)dx$$
where $G(s):=\int_0^s g(\tau)d\tau$ for $s\in \R$. On the other hand, one can use the topological methods, such as the fixed point theory, bifurcation or the Lyapunov-Schmidt reduction. The {\it fixed frequency problem} has been widely studied for the decades, and it is impossible to summarize it here since the related literature is huge.

An important, and of course well known, feature of \eqref{eq:NLSE} is the  conservation of mass:
$$\|\Phi(\cdot,t)\|_{L^2(\R^N)}=\|\Phi(\cdot,0)\|_{L^2(\R^N)}, \quad J[\Phi(\cdot,t)]=J[\Phi(\cdot,0)]\;\hbox{for any}\;t\in \R,$$
where $J$ is the energy functional associated with \eqref{eq:NLSE} defined by
$$J[u]=\frac{1}{2}\int_{\R^N}|\nabla u|^2dx+\frac{1}{2}\int_{\R^N}V(x)u^2dx-\int_{\R^N}G(u)dx,$$
for any
$$u\in \mathcal{H}:=\left\{u\in H^1(\R^N) | \left|\int_{\R^N}V(x)u^2 dx\right|<\infty\right\}.$$
A natural approach to find the  solutions  of \eqref{eq:NLSES} satisfying the normalization constraint
\begin{equation}\label{eq:norm}
 u\in S_a:=\left\{u\in \mathcal{H}\Big| \int_{\R^N}u^2dx=a\right\}
\end{equation}
consists in finding critical points $u\in \mathcal{H}$ of the energy
$J[u]$ under the constraint \eqref{eq:norm}. Then the parameter $\lambda$ appears as a Lagrange multiplier.

\begin{definition}\lab{def:ground}
We write that $u_0$ is a ground state of \eqref{eq:NLSES} on $S_a$ if it is a solution to \eqref{eq:NLSES}  having minimal energy among all the solutions which belongs to $S_a$, i.e.,
$$dJ\big|_{S_a}[u_0]=0\;\hbox{and}\;J[u_0]=\inf\left\{J[u]\Big| dJ\big|_{S_a}[u]=0, u\in S_a\right\}.$$
\end{definition}

\subsection{Non-potential case: $V(x)\equiv 0$}
For $V(x)\equiv 0$, the problem reduces to
\beq\lab{eq:V-tirvial}
\begin{cases}
-\Delta u+\lambda u=g(u)\quad &\hbox{in}\;\R^N,\\
\int_{\R^N}u^2 dx=a>0.
\end{cases}
\eeq
Consider a very special case $g(u)=\sum_{i=1}^{m}a_i|u|^{\sigma_i}u$ with $a_i>0$ and
$$\begin{cases}0<\sigma_i<\frac{4}{N-2}\quad &\hbox{if}\;N\geq 3,\\
\sigma_i>0\;&\hbox{if}\;N=1,2.  \end{cases}\;\forall i=1,2,\cdots,m.$$
For the simplest case $m=1$, $g(u)$ is homogeneous.
Without loss of generality, we assume that $g(u)=|u|^{p-1}u, p\in (1,2^*-1)$, where $\displaystyle 2^*:=\begin{cases}\infty\quad&N=1,2\\ \frac{2N}{N-2} &N\geq 3 \end{cases}$
is the Sobolev critical exponent. That is,
\beq\lab{eq:V-tirvial2}
\begin{cases}
-\Delta u+\lambda u=|u|^{p-1}u\quad &\hbox{in}\;\R^N,\\
\int_{\R^N}u^2 dx=a>0.
\end{cases}
\eeq
Applying the Pohozaev identity, one can show that \eqref{eq:V-tirvial2} possesses a nontrivial solution $u\in H^1(\R^N)$ only if $\lambda>0$. On the other hand,
the positive normalized solution of \eqref{eq:V-tirvial2} can be completely solved by scaling. Let $U_p$ be the unique positive radial solution to
\beq\lab{eq:def-Up}
  -\Delta u+u=u^{p}\;\hbox{in}\;\R^N;\quad u(x)\rightarrow 0\ \text{ as $|x|\to\infty$;}
\eeq
cf.\ \cite{Kwong1989}.
Setting
$$U_{\lambda,p}(x):=\lambda^{\frac{1}{p-1}}U_p(\sqrt{\lambda}x),$$
one can check that, up to a translation, $U_{\lambda,p}$ is the unique positive radial solution to
$$ -\Delta u+\lambda u=u^{p}\;\hbox{in}\;\R^N;\quad u(x)\rightarrow 0\ \text{ as $|x|\to\infty$;}.$$
A direct computation shows that
\beq\lab{eq:U-scaling}
\|U_{\lambda}(x)\|_{L^2(\R^N)}^{2}=\lambda^{\frac{4-(p-1)N}{2(p-1)}}\|U\|_{L^2(\R^N)}^{2}.
\eeq
So we can see that if $p\neq \bar{p}:=1+\frac{4}{N}$, then there exists a  unique $\displaystyle\lambda_a>0$ such that $\displaystyle \|U_{\lambda_a}(x)\|_{L^2(\R^N)}^{2}=a$.
That is, there exists a  positive normalized solution to \eqref{eq:V-tirvial2} for any $a>0$ whenever $p\neq 1+\frac{4}{N}$ (and it is unique up to a translation). While for the so-called mass critical exponent $p+1=2+\frac{4}{N}$, \eqref{eq:V-tirvial2} has positive normalized solution if and only if $a=\|U\|_{L^2(\R^N)}^{2}$ (with infinitely many solutions    and  $\lambda>0$).


If $m\neq 1$,  then $g(u)$ is not of homogeneous, and the scaling method does not work.  Hence, the existence of normalized solutions becomes nontrivial, and many techniques developed for the {\it fixed frequency problem} can not be applied directly. Therefore, the literature focused on the normalized solutions is far less broad. In a series of works \cite{Stuart1980,Stuart1981,Stuart1989}, Stuart applied a bifurcation approach to study the nonhomogeneous nonlinearities (not necessarily autonomous). However, it requires more restrictive growth conditions on $g$ to guarantee the compactness. Stated in the particular case of \eqref{eq:V-tirvial} with $g(u)=\sum_{i=1}^{m}a_i|u|^{\sigma_i}u$,
 Stuart studied the mass sub-critical case, i.e., $0<\sigma_i<\frac{4}{N},1\leq i\leq m$, he showed that the corresponding functional is bounded from  below  on $S_a$ for any $a>0$ and obtained a critical point proving that the infimum is reached (see \cite{Stuart1981,Stuart1989} for the details).
 In \cite{Shibata2014},  M.Shibata studied the mass sub-critical case for general nonlinearities.  Consider
\beq\lab{def:I}
I[u]=\frac{1}{2}\int_{\R^N}|\nabla u|^2 dx-\int_{\R^N}G(u)dx,
\eeq
and study  the $L^2$-constraint minimization problem
\beq\lab{zzz111} \displaystyle E_a=\inf_{u\in\tilde{S}_a}I[u]; \quad\quad
\hbox{ where  } \;\; \displaystyle
\tilde{S}_a:=\left\{u\in H^1(\R^N)\Big| |u|_2^2=a\right\}.
\eeq

 \bl\lab{lemma:Shibata}(cf.\cite[Theorem 1.1 and Theorem 1.3]{Shibata2014})
Assume that
\begin{itemize}
\item[(G1)] $g\in C(\R), g(0)=0$.
\item[(G2)] $\displaystyle \lim_{s\rightarrow 0}\frac{g(s)}{s}=0$.
\item[(G3)] $\displaystyle \lim_{s\rightarrow \infty}\frac{g(s)}{|s|^{1+\frac{4}{N}}}=0$.
\item[(G4)] There exists $s_0>0$ such that $G(s_0)>0$, where $G(s)=\int_0^s g(\tau)d\tau$ for $s\in \R$.
\end{itemize}
Then there exists some $a_0\geq 0$ such that for any $a>a_0$, there exists a global minimizer with respect to $E_a$. That is ,$\exists U\in \tilde{S}_a$ such that
$$dI\big|_{\tilde{S}_a}[U]=0\;\hbox{and}\;I[U]=\inf\left\{I[u]| dI\big|_{\tilde{S}_a}[u]=0, u\in \tilde{S}_a\right\}.$$
In particular, $E_a<0$ for any $a>a_0$ and there is no global minimizer with respect to $E_a$ if $0<a<a_0$. Furthermore, if
$\displaystyle\liminf_{s\rightarrow 0}\frac{g(s)}{|s|^{1+\frac{4}{N}}}=\infty$, then $a_0=0$. If
$\displaystyle\limsup_{s\rightarrow 0}\frac{g(s)}{|s|^{1+\frac{4}{N}}}<\infty$, then $a_0>0$.
\el\hfill$\Box$

While for the mass super-critical case, for example $g(u)=\sum_{i=1}^{m}a_i|u|^{\sigma_i}u$ with $\frac{4}{N}<\sigma_i<2^*-2,1\leq i\leq m$, the corresponding functional is not bounded below anymore. Thus it is impossible to search for a minimum on $S_a$.  L. Jeanjean \cite{Jeanjean1997} could prove the corresponding functional possesses the mountain pass geometric structure on $S_a$. Then obtained  the normalized solution by a minimax approach and a smart compactness argument. Precisely, under the following $(AR)$ type assumption on the  general nonlinearities $g$:
\beq\lab{eq:AR-condition}
\alpha G(s)\leq g(s)s\leq \beta G(s),\;\hbox{where}\;2+\frac{4}{N}<\alpha\leq \beta<2^*,
\eeq
Jeanjean obtained a normalized solution to \eqref{eq:V-tirvial} for any $a>0$ provided $N\geq 2$.
While for the case of $N=1$, to obtain the same result, Jeanjean need the following  additional assumption
\beq\lab{eq:Jeanjean-N=1}
\tilde{G}'(s)s>(2+\frac{4}{N})\tilde{G}(s),\hbox{where}\;\tilde{G}(s)=g(s)s-2G(s).
\eeq
We also note the   combined nonlinearity case,  we refer the readers to \cite{Soave-1,Soave2020} by N.~Soave.

If $g=g(x, u)$ is not autonomous, then the functional is not invariance under transformation anymore. One needs further assumptions to guarantee the compactness.  In \cite{Stuart1982,Stuart1983}, C.A.Stuart studied the special case of $g(x, u)=q(x)|u|^{\sigma} u$, where $q(x)$ is not necessary radial and $\sigma\in (0,\frac{4}{N})$. If $q(x)\rightarrow 0$ as $|x|\rightarrow \infty$ and $q(x)\geq A(1+|x|)^{-t}$ for some $A>0, t\in (0, 2-\frac{N\sigma}{2})$, C.A.Stuart proved that $\lambda=0$ is a bifurcation point for
\beq\lab{eq:Stuart}
-\Delta u+\lambda u=q(x)|u|^{\sigma} u, \quad u\in H^1(\R^N).
\eeq
 And C.A.Stuart raised a problem that whether $\lambda=0$ is still a bifurcation point if $q(x)\not\rightarrow 0$ as $|x|\rightarrow \infty$. Zhu and Zhou in \cite{ZhuZhou1988} solved it under the assumption
\beq\lab{eq:ZhuZhou}
\lim_{u\rightarrow 0}\frac{g(x, u)}{u}=0\;\hbox{uniformly in}\;x\in \R^N
\eeq
and some other assumptions on $g(x, u)$.  We note that $\lambda=0$ is a bifurcation point  implies the existence of normalized solution for small $a>0$.

Recently, Chen and Tang \cite{Chen-Tang2020} study the case of $g(x, u)=q(x)f(u)$ for a general nonlinearity $f(u)$. Under a sequence of technical assumptions, they establish the existence of normalized solutions for both mass super-critical case (see \cite[Theorem 1.1]{Chen-Tang2020}) and mass sub-critical case (see \cite[Theorem 1.4]{Chen-Tang2020}).

\subsection{$V(x)\not\equiv 0$, with potential case}
Observing that in practice physical background, many problems inevitably involve the potentials.
If $V(x)\not\equiv 0$, noting that $-V(x)u+g(u)$ is not  homogeneous, not autonomous. And furthermore,
$$\lim_{u\rightarrow 0}\frac{-V(x)u+g(u)}{u}=-V(x)+\lim_{u\rightarrow 0}\frac{g(x, u)}{u},$$
a value which depends on $x\in \R^N$. So the techniques developed in a sequences of literature mentioned above can not applied directly.
Hence, besides the importance in the applications, not negligible
reasons of our interest for such problems are its stimulating and challenging difficulties.

Recently, in \cite[Section 3]{PellacciPistoiaVairaVerzini}, Pellacci et al. apply Lyapunov-Schmidt reduction approach to study the elliptic equation
\beq\lab{eq:Pellacci}
\begin{cases}
-\Delta v+(V(x)+\lambda)v=v^p\;\hbox{in}\;\R^N,\\
v>0,\int_{\R^N}v^2dx=a.
\end{cases}
\eeq
By putting $\varepsilon:=\lambda^{-\frac{1}{2}}$ and $u=\varepsilon^{\frac{2}{p-1}}v$, \eqref{eq:Pellacci} is reduced to
\beq\lab{eq:Pellacci2}
\begin{cases}
-\varepsilon^2\Delta u+(\varepsilon^2V(x)+1)v=u^p\;\hbox{in}\;\R^N,\\
u>0,\varepsilon^{-\frac{4}{p-1}}\int_{\R^N}u^2dx=a.
\end{cases}
\eeq
Suppose $\xi_0$ is a non-degenerate critical point of $V(x)$, after a  tedious but standard calculation, for $\varepsilon$ small enough, \cite{PellacciPistoiaVairaVerzini} construct a solution $u_\varepsilon$ to \eqref{eq:Pellacci2} with some concentrating behavior such that
$$\varepsilon^{-\frac{4}{p-1}}\int_{\R^N}u_\varepsilon^2dx\rightarrow \begin{cases}\infty\quad&\hbox{if}\;p<1+\frac{4}{N},\\
0&\hbox{if}\;p>1+\frac{4}{N}
\end{cases}\;\hbox{as}\;\varepsilon\rightarrow 0^+.$$
Consequently, they obtain the existence of normalized solution to \eqref{eq:Pellacci} provided $a$ large enough if $p<1+\frac{4}{N}$ and $a$ small enough if $p>1+\frac{4}{N}$.

\br
Consider a general nonlinearity $g(u)$. If the equation
\beq
\begin{cases}
-\Delta u+u=g(u)\;\hbox{in}\;\R^N,\\
u>0, u(x)\rightarrow 0\;\hbox{as}\;|x|\rightarrow \infty
\end{cases}
\eeq
 has an unique solution $U$. And the non-degenerate assumption holds:
 $$Ker(-\Delta +I-g'(U))=span\left\{\frac{\partial U}{\partial x_i}, i=1,2,\cdots,N\right\}.$$
 One can also apply the Lyapunov-Schmidt reduction approach to obtain the existence of normalized solution for $a$ large enough or small enough.
However, for the technique reason, one can not expect to establish the existence result for a large range of $a>0$ by this kind approach.
\er

Consider the mass sub-critical case,
it is well known that,the corresponding functional on the constraint set $S_a$, is bounded below. Then it is nature to search for the minimizer as a critical point of $J$. One of the most essential difficulties in using constrained variational method to study fixed mass problems is to deal with compactness. A key step is to establish the sub-additive inequality that $C_{a+b}<C_{a}+C_{b}$, where
\beq\lab{eq:def-Ca}
C_a:=\inf_{u\in S_a}J[u].
\eeq
According to the traditional argument, the potential $V(x)$ will bring doubling complexity. One of the novelty in present is that we shall give a new approach to deduce the sub-additive inequality, which can greatly simplify the discussion.

\br
\begin{itemize}
\item[(1)]
After finishing the first draft of the paper, we were informed that recently Norihisa Ikoma and Yasuhito Miyamoto in \cite{Ikoma-Miyamoto2020} also study this problem under suitable assumptions. We note that the assumptions they asked for intersected with those we need, but do not contain each other. To make it clear, we shall give a remark after our Theorem \ref{Main-th2} (see Remark \ref{remark:IM}).
\item[(2)]
There are also some progresses on the mass super-critical case, see T.Bartsch et al.\cite{Bartsch2021}. However, they only consider the case of $g(u)=|u|^{p-2}u$.
\end{itemize}
\er

The paper is organized as follows. In the next section we state and discuss our results, in particular we compare them with existing results obtained recently. We shall collect some preliminaries in section \ref{sec:prelim} and  establish some compactness results in section \ref{sec:compactness}. Theorem~\ref{Main-th1} and Theorem~\ref{Main-th2} will be proved in section \ref{sec:proofs}.

\s{Statement of results}
\renewcommand{\theequation}{2.\arabic{equation}}
\renewcommand{\theremark}{2.\arabic{remark}}
\renewcommand{\thedefinition}{2.\arabic{definition}}
The basic idea  is still to find a minimizer of $J$ constrained on $S_a$. However, the procedure is much more complex comparing with the case of $V(x)\equiv 0$. The assumptions we need  are quite  different from those  adopted in \cite{Ikoma-Miyamoto2020}.

\vskip0.12in
Suppose that the nonlinearity $g$ satisfies $(G1)$-$(G4)$. For the potential $V(x)$, we assume that $V(x)\in C(\R^N)$ satisfies
\begin{itemize}
\item[$(V_1)$] $\displaystyle \lim_{|x|\rightarrow \infty}V(x)=\sup_{x\in \R^N}V(x)=:V_\infty\in (0,+\infty]$.
\item[$(V_2)$] $\displaystyle V(0)=\min_{x\in \R^N}V(x)=c_\ell>-\infty$.
\end{itemize}
Firstly, we consider a simple case that $V(x)$ is coercive, i.e., $V_\infty=\infty$. In such a case, the coercive assumption enables us to establish the compactness result. Here comes our first main result.

\bt\lab{Main-th1}
Suppose $(G1)$-$(G4)$ and assume that $V(x)\in C(\R^N)$ satisfies $(V_1)$ with $V_\infty=\infty$ and $(V_2)$.
Then for any $a>0$, there exists a ground state solution  to \eqref{eq:NLSES} with the constraint $u_0\in S_a$, i.e.,
$$dJ\big|_{S_a}[u_0]=0\;\hbox{and}\;J[u_0]=\inf\left\{J[u]\Big| dJ\big|_{S_a}[u]=0, u\in S_a\right\}.$$
\et

Secondly, we study the case of $V_\infty<\infty$. The problem will become very complicated and difficult to handle and therefore the concentration compactness argument (see \cite{Lions1984,Lions1984-2}) have to be involved.  A key step is to prove sub-additive inequality that $ C_{a+b}<C_a+C_b $
 (See  \eqref{eq:def-Ca} ). It is worth mentioning  that for the case of $V(x)\equiv 0$, one can choose $u\in S_a\cap C_c^\infty$ and $v\in S_b\cap C_c^\infty$ such that $supp~u \cap suu~v=\emptyset$, then it is easy to deduce that $E_{a+b}\leq E_a+E_b$ (see \eqref{zzz111}). For the special case that $g(u)=|u|^{p-1}u$ with $1<p<1+\frac{4}{N}$, a direct computation via scaling show that
 $$E_a=\left(\frac{a}{a+b}\right)^{\gamma}E_{a+b}, \forall a,b\in \R^+,$$
 where $\gamma:=\frac{2(p+1)-N(p-1)}{4-(p-1)N}>1$.
 Then by the convex of $q(x):=x^{\gamma}+(1-x)^{\gamma}, x\in [0,1]$, one can see that
 $$p(x)>p(0)=p(1)=1, \forall x\in (0,1),$$
 which implies that strict binding inequality $E_{a+b}<E_a+E_b$ for any $a,b>0$.
 However, when $g$ is not homogenous, the scaling method does not work. For the general case about $g$,  let $u\in S_a$, one can see that $u_\lambda(x):=u(\lambda^{-\frac{1}{N}}x)\in S_{\lambda a}$ and
$$J[u_\lambda]=\lambda \left(\frac{\lambda^{-\frac{2}{N}}}{2}\int_{\R^N}|\nabla u|^2dx -\int_{\R^N}G(u)dx\right).$$
From which one can see that $E_{\lambda a}\leq \lambda E_a$ for any $a>0,\lambda\geq 1.$ Combining with $E_a<0$, one can also finally obtain the strictly monotonicity of $C_a$ and the sub-additive inequality.
Now, let us consider the case of  $V(x)\not\equiv 0$,
\begin{align*}
J[u_\lambda]=\lambda \left(\frac{\lambda^{-\frac{2}{N}}}{2}\int_{\R^N}|\nabla u|^2dx
+\frac{1}{2}\int_{\R^N}V(\lambda^{\frac{1}{N}}x)u^2 dx-\int_{\R^N}G(u)dx\right).
\end{align*}
one can see that $\frac{1}{2}\int_{\R^N}V(\lambda^{\frac{1}{N}}x)u^2 dx$ may increase by $\lambda>1$ while $\frac{\lambda^{-\frac{2}{N}}}{2}\int_{\R^N}|\nabla u|^2dx$ definitely decreases by $\lambda>1$.
This kind of competition will bring great difficulties  for solving the  problem. So the method mentioned above can not  be applied directly.
A natural way is to construct some suitable $\phi\in S_{a+b}$ such that
$J[\phi]<C_a+C_b$. Let $u\in S_a, v\in S_b$ such that $u,v$ nearly attain $C_a$ and $C_b$, then stick u and v together in a proper way. Under some strong assumptions, $J[\phi]<C_a+C_b$ is
expected. However, the process will be very long and tedious.
In present paper, we shall develop a new method basing on iteration to establish the strict binding inequality.
We shall need the following additional $(AR)$-type assumption on the nonlinearity, which has been widely used in studying the fixed frequency problems:
\begin{itemize}
\item[$(G_5)$] There exists some $\alpha>2$ such that $\displaystyle g(s)s\geq \alpha G(s)$ for $s\geq 0$.
\end{itemize}

\br\lab{remark:V3}
Without loss of generality, we may assume that $V_\infty=0$. If not, we may  replace  $(V(x),\lambda)$ by $(\tilde{V}(x),\tilde{\lambda}):=(V(x)-V_\infty,\lambda+V_\infty)$.
\er

\bt\lab{Main-th2}
Suppose $(G1)$-$(G5)$.
Assume that $V(x)\in C^1(\R^N)$ satisfies $(V_1)$,$(V_2)$ and
$$\lim_{|x|\rightarrow \infty}\langle \nabla V(x), x\rangle=0.$$
There exists some $a_0\geq 0$ such that $C_a<0$ for $a>a_0$ and it is attained. While $C_a\equiv 0$ for $0<a<a_0$ and it is not attained. That is,
for any $a>a_0$, there exists
a ground state solution $(u_0,\lambda)$ to \eqref{eq:NLSES} with $u_0\in S_a$ and  $\lambda>-V_\infty$.
\et

\br\lab{remark:zy2}
Let $\alpha_0$ be the sharp number given by \cite[Theorem 1.1]{Shibata2014}, i.e., $E_a<0$ for $a>\alpha_0$ and $E_a\equiv 0$ for $0<a<\alpha_0$. Let $a_0$ be the sharp number given by Theorem \ref{Main-th2} above. Then we have the relation
$a_0\leq \alpha_0.$
That is, $C_a<0$ may happen even $E_a=0$ (It is known that $E_a<0\Rightarrow C_a<E_a<0$).
\er

\bt\lab{Main-th3}
Under the assumptions of Theorem \ref{Main-th2} we suppose further that one of the following holds:
\begin{itemize}
\item[(i)] $\displaystyle \liminf_{t\rightarrow 0} \frac{G(t)}{t^{2+\frac{4}{N}}}=\infty$;
\item[(ii)] There exists an $s_0>0$ such that $g(s)\geq 0$ in $[0,s_0]$ and, in addition if $N \geq 3$,
  $$\inf_{|\varphi|_2^2=1}(|\nabla \varphi|^2+V(x)\varphi^2)dx<0.$$
\end{itemize}
Then $a_0=0$ holds,  that is, for any $a>0$, there exists
a ground state solution $(u_0,\lambda)$ to \eqref{eq:NLSES} with $u_0\in S_a$.
\et

\vskip0.3in
\br\lab{remark:IM} \quad
\begin{itemize}
\item[(1)] For $N\leq 4$, the authors in \cite{Ikoma-Miyamoto2020} need a further assumption on $g$ that:
    \begin{itemize}
    \item [(G6)] $g(s)$ is locally H\"older continuous in  $\R,$  with exponent $\nu\in (0,1)$; $g(s)>0$ for $s>0$ and there exists $\delta_1>0$ such that $g(s)/s$ is nondecreasing in $(0,\delta_1)$.
    \end{itemize}
    However, in our Theorem \ref{Main-th3}, $g(s)\equiv 0$ for $s\in [0,\delta_1]$ with some suitable $\delta_1>0$ is allowed.

\item[(2)]For $N\geq 5$, besides the assumption $(G6)$ above, the authors in \cite{Ikoma-Miyamoto2020} suppose further the following $(G7)$ or $(V_3)$:
    \begin{itemize}
    \item  [(G7)]    $\displaystyle \lim_{s\rightarrow 0} g(s)/|s|^{1+\frac{2}{N-2}}>0$;
    \item [$(V_3)$]  $V\in  W^{1,\infty}(\R^N)$ and $\nabla V(x)\cdot x\leq \frac{(N-2)^2}{2|x|^2}$ for a.e. $x\in \R^N\backslash\{0\}$.
    \end{itemize}
    Then they can prove the existence of normalized solution for $a$ suitable large  $a>a_0$ (see \cite[Theorem A]{Ikoma-Miyamoto2020}).
    While  in the present paper, we do not need $(G6)$.  We also note that $(V_3)$ above is  not necessary for us.

\item[(3)]In \cite[Theorem B-(i)]{Ikoma-Miyamoto2020}, under the assumptions $g(s)\geq 0$ in $[0,s_0]$ and the following $(V_4)$ holds:
\begin{itemize}
    \item  [$(V_4)$] $\;\inf_{\|\varphi\|_{L^2(\R^N)}=1} \int_{\R^N}(|\nabla \varphi|^2+V(x)\varphi^2)dx<0.$
    \end{itemize}
    Then they can prove  the existence of normalized solution for all $a>0$.

\item[(4)]In our opinion, the number $C_a$, which is given by \eqref{eq:def-Ca}, its  negativity depends on $V(x)$ and $g(s)$.  In fact, for $u\in S_a$, we consider that
    \begin{align*}
    &\quad\quad J[u_t]:=J[t^{\frac{N}{2}}u(tx)]=\\
    &\frac{1}{2}t^2 \int_{\R^N} |\nabla u|^2 dx+\frac{1}{2}\int_{\R^N}V(\frac{x}{t})|u|^2dx
-t^{-N}\int_{\R^N}G(t^{\frac{N}{2}}u)dx.
\end{align*}
Noting that the signs of $\frac{1}{2}\int_{\R^N}V(\frac{x}{t})|u|^2dx$ is negative, therefore,
\\
either
\beq\lab{zzz222}
 \frac{1}{2}t^2 \int_{\R^N} |\nabla u|^2 dx+\frac{1}{2}\int_{\R^N}V(\frac{x}{t})|u|^2dx<0\;\hbox{with}\;\int_{\R^N}G(t^{\frac{N}{2}}u)dx\geq 0
 \eeq
or
\beq\lab{zzz333}\frac{1}{2}t^2 \int_{\R^N} |\nabla u|^2 dx-t^{-N}\int_{\R^N}G(t^{\frac{N}{2}}u)dx<0,\eeq
they can imply that $C_a<0$. Hence, we give Theorem \ref{Main-th3} above. Indeed, the assumptions $g(s)\geq 0$ in $[0,s_0]$ and $(V_4)$ are used to guarantee  \eqref{zzz222}  holds for $t$ small  in \cite{Ikoma-Miyamoto2020}. One can also use a condition that $\displaystyle \lim_{s\rightarrow 0} g(s)/|s|^{1+\frac{4}{N}}=\infty$ to guarantee the inequality  \eqref{zzz333}  happens provided $t>0$ small enough, which is considered in \cite{Shibata2014}.

\item[(5)]In addition, we need to point out that the conditions $(V_3)$ and $(V_4)$ are mutually exclusive,  which means that for many cases, they can not obtain the result $a_0=0$. Indeed, take $V(x)$, a radial function, as an example. If $(V_3)$ holds, one have that $V'(r)\leq \frac{(N-2)^2}{2r^3},$ combing with $\displaystyle \lim_{|x|\rightarrow \infty}V(x)=0$, one can prove that
    $$\int_s^\infty V'(r)dr\leq \int_s^\infty \frac{(N-2)^2}{2r^3}dr\Rightarrow V(x)\geq -\frac{(N-2)^2}{4|x|^2}.$$
    Then, by Hardy inequality,
    $$\int_{\R^N}(|\nabla \varphi|^2+V(x)\varphi^2)dx\geq \int_{\R^N}(|\nabla \varphi|^2-\frac{(N-2)^2}{4|x|^2}\varphi^2)dx\geq 0,$$
    which means that hypothesis $(V_4)$ is impossible provided $(V_3)$. However, such a situation will not happen in our present paper.
\end{itemize}
\er

\br\lab{remark:no}
One  novelty in the current work is presenting a new approach to establish the sub-additive inequality, which is a crucial step to study the mass-subcritical problem (and we believe it can be applied to study the non-potential case, even for the system case). The assumptions on the potential $V(x)$ are quite relax in the present paper.
\er

\s{Preliminaries }\lab{sec:prelim}
\renewcommand{\theequation}{3.\arabic{equation}}
\renewcommand{\theremark}{3.\arabic{remark}}
\renewcommand{\thedefinition}{3.\arabic{definition}}
We write $|u|_p$ for the $L^p$-norm.
Since the frequency $\lambda$ is unknown, we can not use the usual Nehari manifold. However, combing the usual Nehari manifold and the Pohozaev manifold, we can obtain the following result.

\bl\lab{lemma:Pohozaev}
Assume that $V(x)\in C^1(\R^N)$ and $u\in H^1(\R^N)$ is a solution to \eqref{eq:NLSES} with constraint \eqref{eq:norm}, then
$u\in S_a\cap \mathcal{P}$, where
\begin{align*}
\mathcal{P}:=&\left\{u\in H^1(\R^N) \Big| \int_{\R^N}|\nabla u|^2dx-\frac{1}{2}\int_{\R^N}\langle \nabla V(x), x\rangle u^2 dx\right.\\
&\quad\left.+N\int_{\R^N}\left[G(u)-\frac{1}{2}g(u)u\right]dx=0\right\}.
\end{align*}
\el

\bp
Let $u$ be a solution to \eqref{eq:NLSES}, then we have
\beq\lab{eq:Nehari}
\int_{\R^N}|\nabla u|^2dx+\int_{\R^N}(V(x)+\lambda)u^2dx -\int_{\R^N}g(u)udx=0
\eeq
and
\begin{align}\lab{eq:Pohozaev}
&(N-2)\int_{\R^N}|\nabla u|^2dx+N\int_{\R^N}(V(x)+\lambda)u^2dx  \nonumber \\
&-2N\int_{\R^N}G(u)dx+\int_{\R^N}\langle \nabla V(x), x\rangle u^2 dx=0.
\end{align}
Eliminate the unknown parameter $\lambda$, we obtain that
$$\int_{\R^N}|\nabla u|^2dx-\frac{1}{2}\int_{\R^N}\langle \nabla V(x), x\rangle u^2 dx+N\int_{\R^N}\left[G(u)-\frac{1}{2}g(u)u\right]dx=0.$$
\ep

\bl\lab{lemma:J-bounded-blow}
Suppose $(G1)$-$(G4)$   and $V(x)\geq c_\ell$ in $\R^N$.  Then  $ J[u]$ is bounded below on $S_a$, i.e.,
$$C_a:=\inf_{u\in S_a}J[u]>-\infty.$$
Furthermore, if $V_\infty<\infty$, it holds that $C_a\leq E_a+\frac{1}{2}V_\infty a$. In particular,  $C_a<E_a+\frac{1}{2}V_\infty a$ if $E_a$ is attained.
\el
\bp
Noting that $J[u]=I[u]+\frac{1}{2}\int_{\R^N}V(x)u^2dx$, by $V(x)\geq c_\ell$, we have that
$$C_a:=\inf_{u\in S_a}J[u]\geq \inf_{u\in S_a}I[u]+\frac{c_\ell}{2}a\geq \inf_{u\in \tilde{S}_a}I[u]+\frac{c_\ell}{2}a=E_a+\frac{c_\ell}{2}a>-\infty.$$
If $V_\infty=\infty$, then it is trivial that $C_a<E_a+\frac{1}{2}V_\infty a=\infty$. If $V_\infty<\infty$, then we have that $\mathcal{H}=H^1(\R^N)$ and $S_a=\tilde{S}_a$.
For any $\varepsilon>0$, we can take $u\in S_a$ such that $I[u]<E_a+\varepsilon$. Observe that for any $R>0$, $u(\cdot-R)\in S_a$ and $I[u(\cdot)]=I(u(\cdot-R))$.
So
$$C_{a}\leq \lim_{R\rightarrow \infty}J[u(\cdot-R)]$$
$$=\lim_{R\rightarrow \infty}\Big[I(u(\cdot-R))+\frac{1}{2}\int_{\R^N}V(x)u(x-R)^2dx\Big]\leq E_a+\frac{1}{2}V_\infty a+\varepsilon.$$
By the arbitrariness of $\varepsilon$, we obtain that
$C_a\leq E_a+\frac{1}{2}V_\infty a$. If $E_a$ is attained, let $U$ be a minimizer of $I$ on $S_a$.
By $V(x)\leq V_\infty$ but $V(x)\not\equiv V_\infty$, we have that
\begin{align*}
C_a\leq &J[U]=I[U]+\frac{1}{2}\int_{\R^N}V(x)Udx\\
=&E_a+\frac{1}{2}\int_{\R^N}V(x)U^{2}dx\\
<&E_a+\frac{1}{2}\int_{\R^N}V_\infty U^{2}dx=E_a+\frac{1}{2}V_\infty a.
\end{align*}
\ep

\bl\lab{lemma:bounded-bounded}
Suppose $(G1)$-$(G3)$ and
$V(x)\geq c_\ell$ in $\R^N$.
Then any $L^2$-bounded sequence $\{u_n\}$ with $J[u_n]<\infty$ is also bounded in $H^1(\R^N)$.
\el
\bp
Let $\{u_n\}$ be a $L^2$-bounded sequence, i.e., $|u_n|_2^2\leq M_1$. Suppose that
$J[u_n]\leq M_2,\;\forall \;n$.
That is,
\beq\lab{eq:Nehari-identity}
\frac{1}{2}|\nabla u_n|_2^2dx+\frac{1}{2}\int_{\R^N}V(x)u_n^2dx-\int_{\R^N}G(u_n)dx\leq M_2.
\eeq
 We shall prove that $u_n$ is also bounded in $H^1(\R^N)$.
By the assumptions $(G1)$-$(G3)$, we have that for any $\varepsilon>0$, there exists some $C_\varepsilon>0$  such that
\beq\lab{eq:Young-inequality}
G(y)\leq C_\varepsilon y^2+\varepsilon y^{2+\frac{4}{N}}, \forall y\geq 0.
\eeq
So we have
\beq\lab{eq:Young-inequality-2}
\int_{\R^N}G(u_n)dx\leq C_\varepsilon |u_n|_2^2+\varepsilon |u_n|_{2+\frac{4}{N}}^{2+\frac{4}{N}}.
\eeq
By Gagliardo-Nirenberg inequality, there exists some $C(N)>0$ such that we have
\beq\lab{eq:Gagliardo-Nirenberg}
|u|_{2+\frac{4}{N}}^{2+\frac{4}{N}}\leq C(N)|\nabla u|_2^2 |u|_{2}^{\frac{4}{N}},\;\forall u\in H^1(\R^N).
\eeq
Then by \eqref{eq:Nehari-identity}, \eqref{eq:Young-inequality-2}, \eqref{eq:Gagliardo-Nirenberg}, we
obtain that
\begin{align*}
\frac{1}{2}|\nabla u_n|_2^2\leq &M_2-\frac{1}{2}\int_{\R^N}V(x)u_n^2dx+\int_{\R^N}G(u_n)dx\\
\leq&M_2-\frac{1}{2}c_\ell M_1+\left[C_\varepsilon |u_n|_2^2+\varepsilon |u_n|_{2+\frac{4}{N}}^{2+\frac{4}{N}}\right]\\
\leq&M_2-\frac{1}{2}c_\ell M_1+C_\varepsilon M_1+\varepsilon C(N)|\nabla u_n|_2^2 |u_n|_{2}^{\frac{4}{N}}\\
\leq &M_2-\frac{1}{2}c_\ell M_1+C_\varepsilon M_1+\varepsilon  C(N) M_{1}^{\frac{2}{N}} |\nabla u_n|_2^2.
\end{align*}
So we can take some suitable small $\varepsilon>0$ such that
$$\varepsilon  C(N) M_{1}^{\frac{2}{N}}<\frac{1}{4},$$
 it follows that
\beq\lab{eq:bounded}
|\nabla u_n|_2^2\leq 4\left[M_2-\frac{1}{2}c_\ell M_1+C_\varepsilon M_1\right]<\infty.
\eeq
\ep

\bl\lab{lemma:ls}
Suppose $(G2)$ and $(G3)$. Assume that $V(x)\in C^1(\R^N)$ and
$V(x)\geq c_\ell$ in $\R^N$.
Then for any $u\in S_a$ with $J[u]<\frac{1}{2}V_\infty a$, there exists some $t_0>0$ such that $u_{t_0}\in S_a\cap \mathcal{P}$ and
$$J[u_{t_0}]=\min_{t>0}J[u_t],$$
where $u_t(x)$ is defined by
\beq\lab{eq:ut-invariant}
u_t(x):=t^{\frac{N}{2}}u(tx), t>0.
\eeq
\el
\bp
By a direct computation, we have
$$
\begin{cases}
\int_{\R^N} |\nabla u_t|^2 dx=t^2 \int_{\R^N} |\nabla u|^2 dx,\\
\int_{\R^N}V(x)|u_t|^2dx=\int_{\R^N}V(\frac{x}{t})|u|^2dx,\\
\int_{\R^N}G(u_t)dx=t^{-N}\int_{\R^N}G(t^{\frac{N}{2}}u)dx,
\end{cases}
$$
and then
\begin{align*}
J[u_t]=&\frac{1}{2}t^2 \int_{\R^N} |\nabla u|^2 dx+\frac{1}{2}\int_{\R^N}V(\frac{x}{t})|u|^2dx
-t^{-N}\int_{\R^N}G(t^{\frac{N}{2}}u)dx.
\end{align*}
Under the assumptions of $(G2)$ and $(G3)$,we have that
$$G(t^{\frac{N}{2}}u)=o(1)(t^{\frac{N}{2}}u)^{2+\frac{4}{N}}=o(1)t^{N+2}u^{2+\frac{4}{N}}\;\hbox{as}\;t\rightarrow \infty$$
and
$$G(t^{\frac{N}{2}}u)=o(1)(t^{\frac{N}{2}}u)^{2}=o(1)t^{N}u^{2}\;\hbox{as}\;t\rightarrow 0.$$
Hence,
$$\lim_{t\rightarrow 0} J[u_t]=\frac{1}{2}V_\infty a\;\hbox{and}\;\lim_{t\rightarrow \infty} J[u_t]=+\infty.$$
Noting that $\displaystyle J[u]<\min\{\lim_{t\rightarrow 0} J[u_t],\lim_{t\rightarrow \infty} J[u_t]\}$, we see that there exists some $t_0\in (0,+\infty)$ such that
$$J(u_{t_0})=\min_{t>0}J[u_t]\leq J[u].$$
Then we have  $\frac{d}{dt}J[u_t]\Big|_{t=t_0}=0$, i.e.,$u_{t_0}\in \mathcal{P}$. Hence, $u_{t_0}\in S_a\cap\mathcal{P}$.
\ep

\bc\lab{cro:bounded-away}
Under the assumptions of Lemma \ref{lemma:ls},
for any $a>0$, and $\varepsilon>0$, there exists some $\delta(a,\varepsilon),\tilde{\delta}_{a,\varepsilon}>0$ such that
$$|\nabla u|_2^2\geq \delta(a,\varepsilon)\;\hbox{and}\;\int_{\R^N}G(u)dx\geq \tilde{\delta}_{a,\varepsilon},\;  \forall u\in S_a\cap \mathcal{P}\;\hbox{with}\;J[u]<\frac{1}{2}V_\infty a-\varepsilon.$$
\ec
\bp
We argue by contradiction and assume that there exists $\{u_n\}\subset S_a$ with $|\nabla u_n|_2^2=1$, and  a sequence $t_n\rightarrow 0$ such that $(u_n)_{t_n}\in \mathcal{P}$ and
$$\lim_{n\rightarrow \infty}J[(u_n)_{t_n}]\leq \frac{1}{2}V_\infty a-\varepsilon<\frac{1}{2}V_\infty a.$$
Then by the proof of Lemma \ref{lemma:ls}, we have that as $t\rightarrow 0$,
$$
\begin{cases}
\int_{\R^N} |\nabla u_t|^2 dx=t^2 \int_{\R^N} |\nabla u|^2 dx=t^2=o(1),\\
\int_{\R^N}V(x)|u_t|^2dx=\int_{\R^N}V(\frac{x}{t})|u|^2dx=V_\infty a+o(1),\\
\int_{\R^N}G(u_t)dx=t^{-N}\int_{\R^N}G(t^{\frac{N}{2}}u)dx=\int_{\R^N}o(1)u^2dx=o(1).
\end{cases}
$$
Hence,
$$\lim_{n\rightarrow \infty}J[(u_n)_{t_n}]=\frac{1}{2}V_\infty a >\frac{1}{2}V_\infty a-\varepsilon,$$
a contradiction.

By a similar argument, we can assume that $\{u_n\}\subset S_a$ with $\int_{\R^N}G(u_n)dx=1$, and  a sequence $t_n\rightarrow 0$ such that $(u_n)_{t_n}\in \mathcal{P}$. And then we can also obtain that there exists some $\tilde{\delta}_{a,\varepsilon}>0$ such that
$$\int_{\R^N}G(u)dx\geq \tilde{\delta}_{a,\varepsilon}\; \forall u\in S_a\cap \mathcal{P}\;\hbox{with}\;J[u]<\frac{1}{2}V_\infty a-\varepsilon.$$
\ep

\s{Compactness results}\lab{sec:compactness}
\renewcommand{\theequation}{4.\arabic{equation}}
\renewcommand{\theremark}{4.\arabic{remark}}
\renewcommand{\thedefinition}{4.\arabic{definition}}
\subsection{The coercive case: $V_\infty=\infty$}\lab{sec:infinity}
In this section, we consider $N\geq 1$ and study the coercive case that $V_\infty=\infty$.

\bt\lab{thm:compact-infinity}
Suppose $(G1)$-$(G4)$.
Assume $V(x)\in C^1(\R^N)$ satisfies $(V_1)$ and $(V_2)$ with $V_\infty=\infty$. Then any minimizing sequence $\{u_n\}\subset S_a$ with $J[u_n]\rightarrow C_a$, possesses a convergent subsequence in $H^1(\R^N)$.
\et
\bp
By Lemma \ref{lemma:J-bounded-blow}, $J$ is bounded below on $S_a$.
Let $\{u_n\}\subset S_a$ be a minimizing sequence of $J$ respect to $C_a$. By Lemma \ref{lemma:bounded-bounded}, we have that $\{u_n\}$ is bounded in $H^1(\R^N)$. So we may assume that $u_n\rightharpoonup u$ in $\mathcal{H}$ along a subsequence, and also holds in $H^1(\R^N)$.

\vskip 0.1in\noindent
{\bf Claim 1:} $|u|_2^2=a$ and thus $u_n\rightarrow u$ in $L^2(\R^N)$.\\
If not, we put $v_n=u_n-u$, then $v_n\rightharpoonup 0$ but $v_n\not\rightarrow 0$ in $H^1(\R^N)$.
So that
$$\liminf_{n\rightarrow \infty}|v_n|_2^2=:\delta>0.$$
Up to a subsequence, we may assume that $v_n\rightarrow 0$ in $L_{loc}^{2}(\R^N)$. Then we have that
$$\int_{\R^N}V(x)|v_n|^2 dx\rightarrow \infty,$$
which implies that
$$\int_{\R^N}V(x)|u_n|^2 dx\rightarrow \infty.$$
Recalling that $\{u_n\}$ is bounded in $H^1(\R^N)$, we obtain that
$$C_a=J[u_n]+o(1)=I[u_n]+\frac{1}{2}\int_{\R^N}V(x)|u_n|^2 dx\geq E_a+\int_{\R^N}V(x)|u_n|^2 dx\rightarrow +\infty,$$
a contradiction.
So we obtain that $u_n\rightarrow u$ in $L^2(\R^N)$ and thus $|u|_2^2=a$.

\vskip 0.1in\noindent
{\bf Claim 2:} $u_n\rightarrow u$ in $L^p(\R^N)$ for any $2<p<2^*$.\\
We assume that $|\nabla u_n|_2^2\leq C$ for all $n\in \N$.   By Gagliardo-Nirenberg inequality, we have that
$$|u_n-u|_p^p\leq C(N)|\nabla (u_n-u)|_{2}^{\frac{N(p-2)}{2}} |u_n-u|_{2}^{\frac{2N-(N-2)p}{2}},$$
where $C(N)$ is a positive constant which depends on the dimension $N$.
Hence, the claim follows by the fact of $u_n\rightarrow u$ in $L^2(\R^N)$.
Consequently, under the assumptions $(G1)$-$(G3)$, we have that
$$\int_{\R^N}G(u_n)dx\rightarrow \int_{\R^N}G(u)dx\;\hbox{as}\;n\rightarrow \infty.$$
\vskip 0.1in\noindent
{\bf Claim 3:} $u_n\rightarrow u$ in $H^1(\R^N)$.
If not, we have that
$$\liminf_{n\rightarrow \infty}|\nabla u_n|_2^2>|\nabla u|_2^2.$$
Noting that $u_n\rightharpoonup u$ in $\mathcal{H}$ implies that $$\int_{\R^N}V(x)u_n^2dx=\int_{\R^N}V(x)u^2dx+\int_{\R^N}V(x)(u-u_n)^2dx+o(1).$$
By $u_n\rightarrow u$ in $L_{loc}^{2}(\R^N)$ and $V(x)\in C^1(\R^N)$, we have that
$$\lim_{n\rightarrow \infty} \int_{\R^N}V(x)(u-u_n)^2dx\geq 0,$$
and thus
$$\liminf_{n\rightarrow \infty}\int_{\R^N}V(x)u_n^2dx\geq \int_{\R^N}V(x)u^2dx.$$
So we have
\begin{align*}
C_a=&\lim_{n\rightarrow \infty}\left[\frac{1}{2}|\nabla u_n|_2^2+\frac{1}{2}\int_{\R^N}V(x)u_n^2 dx-\int_{\R^N}G(u_n)\right]\\
>&J[u]\geq C_a,
\end{align*}
a contradiction. Hence, we have
$$\liminf_{n\rightarrow \infty}\int_{\R^N}V(x)u_n^2dx= \int_{\R^N}V(x)u^2dx$$
and
$$\liminf_{n\rightarrow \infty}|\nabla u_n|_2^2= |\nabla u|_2^2.$$
In a word, $u_n\rightarrow u$ both in $H^1(\R^N)$ and in $\mathcal{H}$. Furthermore, $u$ is a minimizer, i.e., $u\in S_a$ and $J[u]=C_a$.
\ep

\vskip0.3in

\subsection{The case of $V_\infty<\infty$}\lab{sec:decreasing}

As mentioned in Remark \ref{remark:V3}, we only need to consider the case of $V_\infty=0$ and  so that $V(x)\leq 0$.

\bl\lab{lemma:decreasing}
Suppose $(G1)$-$(G5)$. Assume that
the potential $V(x)\in C(\R^N)$ satisfies the assumptions $(V_1)$, $(V_2)$ and
$$\lim_{|x|\rightarrow \infty}\langle \nabla V(x),x\rangle=0.$$
Then the following hold.
\begin{itemize}
\item[(1)] $C_a\leq E_a\leq 0$ for $a>0$;  and $C_a<E_a$ if $E_a$ is attained.
\item[(2)] $C_{a+b}\leq C_a+C_b$ for $a,b>0$.   If $C_a$ or $C_b$ is attained,  then $C_{a+b}<C_a+C_b$ for $a, b>0$. In addition,  $C_{\lambda a} < \lambda C_a$ for all $\lambda >1$, if $C_a$ is reached.
\item[(3)] $a\mapsto C_a$ is nonincreasing.
\item[(4)]For sufficiently large $a$, $C_a<0$ holds.
\item[(5)] $a\mapsto C_a$ is continuous.
\end{itemize}
\el
\bp
Observe that $(1)$ and $(4)$ are direct conclusions of \cite[Lemma 2.3]{Shibata2014}.

(3) Now we prove that $C_b\leq C_a$ for any $b>a>0$. Fix any $a>0$ and for any $\varepsilon>0$, we can take some positive functions $u\in S_a\cap C_0^\infty(\R^N)$ and $v\in S_{b-a}\cap  C_0^\infty(\R^N)$ such that
$$J[u]<C_a+\frac{\varepsilon}{2}, \; \;  I[v]\leq E_{b-a}+\frac{\varepsilon}{2}.$$
Since $u$ and $v$ have compact support, we can take $R$ large enough such that
$$\tilde{v}(x):=v(x-R), supp\;u\cap supp\;\tilde{v}=\emptyset.$$
Therefore, $u+\tilde{v}\in S_{b}$. Thus,
$$C_{b}\leq J[u+\tilde{v}]=J[u]+J[\tilde{v}]<J[u]+I[v]\leq C_a+\frac{\varepsilon}{2}+E_{b-a}+\frac{\varepsilon}{2}\leq C_a+\varepsilon,$$
here we use the fact $E_{b-a}\leq 0$.
Then by the arbitrariness of $\varepsilon$, we obtain that $C_b\leq C_a$ for any $b>a>0$.

\vskip0.12in

(5) We only need to prove the following  two claims.

\vskip0.12in
 {\bf Claim 1:} $\displaystyle \lim_{h\rightarrow 0^+}C_{a-h}=C_a$.\\
If $C_a=0$, then by (1) and $(3)$, we have that
$C_{a-h}\equiv 0, \;\forall h\in (0, a)$,
and thus the Claim 1 is trivial. So we consider the case $C_a<0$. For any $u\in S_a$, we put
$u_h(x):=\sqrt{1-\frac{h}{a}}u(x)$ for $h\in (0,a)$. Observe that $u_h\in S_{a-h}$ and
$u_h\rightarrow u$ in $H^1(\R^N)$ as $h\rightarrow 0^+$, we have
$J[u_h]\rightarrow J[u]\;\hbox{as}\;h\rightarrow 0^+$.
Hence,
$$\lim_{h\rightarrow 0^+}C_{a-h}\leq \lim_{h\rightarrow 0^+}J[u_h]=J[u].$$
By the arbitrariness  of $u\in S_a$, we obtain that
$\displaystyle\lim_{h\rightarrow 0^+}C_{a-h}\leq C_a$.
Note that the inequality in the opposite direction is implied by (3), hence the Claim 1 is proved.

\vskip0.12in
{\bf Claim 2:} $\displaystyle \lim_{h\rightarrow 0^+}C_{a+h}=C_a$.\\
We note that the monotonicity property (3) implies the convergence of the  left hand side. It is sufficient to consider the case $h=\frac{1}{n}, n\in \N$. Take $u_n\in S_{a+\frac{1}{n}}$ such that $J[u_n]<C_{a+\frac{1}{n}}+\frac{1}{n}$. Then by Lemma \ref{lemma:bounded-bounded}, we see that $\{u_n\}$ is bounded in $H^1(\R^N)$. Put
$$v_n(x):=\sqrt{\frac{na}{na+1}}u_n(x).$$
We note that $v_n\in S_a$ and $$\|v_n-u_n\|_{H^1(\R^N)}=\left(1-\sqrt{\frac{na}{na+1}}\right)\|u_n\|_{H^1(\R^N)}\rightarrow 0\;\hbox{as}\;n\rightarrow \infty,$$
which implies that
$$J[v_n]=J[u_n]+o(1).$$
Hence,
$$C_a\leq \liminf_{n\rightarrow \infty}J[v_n]=\liminf_{n\rightarrow \infty}[J[u_n]+o(1)]=\lim_{h\rightarrow 0^+}C_{a+h}.$$
The inequality in the opposite direction is also implied by (3). Hence, Claim 2 is proved.

By Claim 1 and Claim 2, we finish the proof of (5).

\vskip0.12in

(2) If $C_a=0$, then by (1) and (3), we have that
$$C_{\lambda a}\leq C_a=0=\lambda C_a\;\hbox{for}\;\lambda>1.$$
So we consider the case $C_a<0$.

Firstly, we prove that
$$C_{\lambda a}\leq\lambda C_a\;\hbox{for any $a>0$ and for $\lambda>1$ closing  to $1$.  }$$
For any $\varepsilon>0$, we can take $u\in S_a\cap \mathcal{P}$ such that
$$J[u]<C_a+\varepsilon.$$
Put $\tilde{u}(\tau, x):=u(\tau^{-\frac{1}{N}}x)$ for $\tau\geq 1$, we have $|\tilde{u}(\tau, x)|_2^2=\tau a$ and
\begin{align*}
J[\tilde{u}(\tau, x)]=&\frac{1}{2}|\nabla \tilde{u}(\tau, x)|_2^2+\frac{1}{2}\int_{\R^N}V(x)\tilde{u}(\tau, x)^2dx-\int_{\R^N}G(\tilde{u}(\tau, x))dx\\
=&\frac{1}{2}\tau^{\frac{N-2}{N}}|\nabla u|_2^2+\frac{\tau}{2}\int_{\R^N}V(\tau^{\frac{1}{N}}x)|u|^2 dx-\tau \int_{\R^N}G(u)dx.
\end{align*}
By a direct computation,
\begin{align*}
\frac{d}{d\tau}J[\tilde{u}(\tau, x)]=&\frac{N-2}{2N}\tau^{-\frac{2}{N}}|\nabla u|_2^2-\int_{\R^N}G(u)dx\\
&+\frac{1}{2}\int_{\R^N}\left[V(\tau^{\frac{1}{N}}x)+\frac{1}{N}
\langle\nabla V(\tau^{\frac{1}{N}}x), \tau^{\frac{1}{N}}x\rangle\right]|u|^2dx.
\end{align*}
Since $u\in \mathcal{P}$, we have
$$\int_{\R^N}|\nabla u|^2dx=\frac{1}{2}\int_{\R^N}\langle \nabla V(x), x\rangle u^2 dx-N\int_{\R^N}\left[G(u)-\frac{1}{2}g(u)u\right]dx.$$
Then
{\allowdisplaybreaks
\begin{align*}
&\frac{d}{d\tau}J[\tilde{u}(\tau, x)]-J[u]=[\frac{N-2}{2N}\tau^{-\frac{2}{N}}-\frac{1}{2}]|\nabla u|_2^2\\
&\quad \quad \quad \quad \quad \quad +\frac{1}{2}\int_{\R^N}\left[V(\tau^{\frac{1}{N}}x)-V(x)+\frac{1}{N}\langle\nabla V(\tau^{\frac{1}{N}}x), \tau^{\frac{1}{N}}x\rangle\right]|u|^2dx\\
=&[\frac{N-2}{2N}\tau^{-\frac{2}{N}}-\frac{1}{2}]\left\{\frac{1}{2}\int_{\R^N}\langle \nabla V(x), x\rangle u^2 dx-N\int_{\R^N}\left[G(u)-\frac{1}{2}g(u)u\right]dx\right\}\\
&\quad \quad \quad \quad \quad \quad +\frac{1}{2}\int_{\R^N}\left[V(\tau^{\frac{1}{N}}x)-V(x)+\frac{1}{N}\langle\nabla V(\tau^{\frac{1}{N}}x), \tau^{\frac{1}{N}}x\rangle\right]|u|^2dx\\
=&\frac{1}{2}\int_{\R^N}\left[V(\tau^{\frac{1}{N}}x)-V(x)+\frac{1}{N}\langle\nabla V(\tau^{\frac{1}{N}}x), \tau^{\frac{1}{N}}x\rangle-\frac{1}{N}\langle\nabla V(x), x\rangle\right]|u|^2dx\\
&+(\frac{N-2}{4N}\tau^{-\frac{2}{N}}-\frac{1}{4}+\frac{1}{2N}) \int_{\R^N}\langle\nabla V(x), x\rangle|u|^2dx\\
&-(\frac{N-2}{2}\tau^{-\frac{N}{2}}-\frac{N}{2}) \int_{\R^N}[G(u)-\frac{1}{2}g(u)u]dx.
\end{align*}}
 By $\displaystyle V(x)\in C^1(\R^N),\lim_{|x|\rightarrow \infty}V(x)=0$ and $\displaystyle \lim_{|x|\rightarrow \infty }\langle \nabla V(x), x\rangle =0$, there exists some $\delta>0$ such that
 $$\sup_{x\in \R^N}\left[V(\tau^{\frac{1}{N}}x)-V(x)+\frac{1}{N}\langle\nabla V(\tau^{\frac{1}{N}}x), \tau^{\frac{1}{N}}x\rangle-\frac{1}{N}\langle\nabla V(x), x\rangle\right]<\frac{\alpha-2}{4a}\tilde{\delta}_{a,\varepsilon}$$
 and
$$\Big(\frac{N-2}{4N}\tau^{-\frac{2}{N}}-\frac{1}{4}+\frac{1}{2N}\Big) \sup_{x\in \R^N}[\langle\nabla V(x), x\rangle]<\frac{\alpha-2}{4a}\tilde{\delta}_{a,\varepsilon}$$
uniformly for any $\tau\in [1,1+\delta]$, where $\varepsilon>0$ is chosen small enough such that $C_a+\varepsilon<0$ and  $\tilde{\delta}_{a,\varepsilon}$ is given by Corollary \ref{cro:bounded-away}.
Then we see that
\begin{align*}
&\frac{d}{d\tau}J[\tilde{u}(\tau, x)]-J[u]\\
<&\frac{1}{2}\frac{\alpha-2}{4a}\tilde{\delta}_{a,\varepsilon_0}\cdot a+\frac{\alpha-2}{4a}\tilde{\delta}_{a,\varepsilon_0}\cdot a+\int_{\R^N}[G(u)-\frac{1}{2}g(u)u]dx\\
\leq&\frac{\alpha-2}{8}\tilde{\delta}_{a,\varepsilon_0}+\frac{\alpha-2}{4}\tilde{\delta}_{a,\varepsilon_0}
-\frac{\alpha-2}{2}\int_{\R^N}G(u)dx\\
<&-\frac{1}{8}\tilde{\delta}_{a,\varepsilon}<0\;\hbox{for any}\;\tau\in [1,1+\delta].
\end{align*}
We note that the $\delta$ can be chosen independent of the choice of $u\in S_a\cap \mathcal{P}$ (certainly it depends on $a$). Indeed, we can firstly take some $\varepsilon_0>0$ small enough such that $C_a+\varepsilon_0<0$. Then by Corollary \ref{cro:bounded-away}, there exists some $\delta_{a,\varepsilon_0},\tilde{\delta}_{a,\varepsilon_0}>0$ such that
$$0< \delta_{a,\varepsilon_0}\leq |\nabla u|_2^2\;\hbox{and}\;0<\tilde{\delta}_{a,\varepsilon_0}\leq\int_{\R^N}G(u)dx$$
for any $u\in S_a\cap \mathcal{P}$ with  $J[u]<C_a+\varepsilon$ and $\varepsilon\in (0,\varepsilon_0]$.
Hence, for any $\lambda\in (1,1+\delta)$, we have that
$$J[\tilde{u}(\lambda, x)]-J[u]=\int_1^\lambda \frac{d}{d\tau}J[\tilde{u}(\tau, x)] d\tau<\int_1^\lambda J[u]d\tau=J[u](\lambda -1),$$
and thus
$$C_{\lambda a}\leq J[\tilde{u}(\lambda, x)]\leq \lambda J[u]\leq \lambda (C_a+\varepsilon).$$
Finally, by the arbitrariness  of $\varepsilon$,  \; we obtain that
$$C_{\lambda a}\leq \lambda C_a, \lambda\in (1,1+\delta).$$
It is necessary to point out that if $C_a$ is attained, one can take $u$ as an minimizer in the argument above, and obtain the strictly inequality
$$C_{\lambda a}< J[\tilde{u}(\lambda, x)]\leq \lambda J[u]=\lambda C_a,, \lambda\in (1,1+\delta).$$

\vskip0.1in

Furthermore, following the proof of Corollary \ref{cro:bounded-away}, since $C_a$ is nonincreasing, if $C_a<0$, for any $b\in (a,+\infty)$, we can find some uniform $\delta>0$ such that
$$C_{\lambda c}\leq \lambda C_c, \forall \lambda\in [1,1+\delta), \forall c\in [a,b],$$
due to the fact $C_b+\varepsilon\leq C_a+\varepsilon<0$. Now, for any $a>0$ with $C_a<0$ and any $\lambda>1$. We take $\delta>0$ such that
$$C_{(1+t) c}\leq (1+t) C_c, \forall t\in [0,\delta), \forall c\in [a,\lambda a].$$
Then we can choose $t_0\in (0,\delta)$ and $m\in \N$ such that
$$(1+t_0)^m\leq \lambda <(1+t_0)^{m+1},$$
and then
\begin{align*}
C_{\lambda a}=&C_{(1+t_0)\frac{\lambda}{1+t_0}a}\leq (1+t_0)C_{\frac{\lambda}{1+t_0}a}\\
\leq&(1+t_0)^2C_{\frac{\lambda}{(1+t_0)^2}a}\\
\leq&\cdots\\
\leq&(1+t_0)^mC_{\frac{\lambda}{(1+t_0)^m}a}\begin{cases} =\lambda C_a\;&\hbox{if}\;\lambda=(1+t_0)^m,\\
\leq(1+t_0)^m\frac{\lambda}{(1+t_0)^m}C_{a}=\lambda C_a&\hbox{if}\;\lambda>(1+t_0)^m.
\end{cases}
\end{align*}
Furthermore, if $C_a$ is attained, then the final step above is strictly ``$<$", and we obtain that $C_{\lambda a}<\lambda C_a$ for any $\lambda>1$.

Hence, for any $0<b\leq a$, we have that
$$C_{a+b}=C_{\frac{a+b}{a}a}\leq \frac{a+b}{a}C_a=C_a+\frac{b}{a}C_a=C_a+\frac{b}{a}C_{\frac{a}{b}b}\leq C_a+\frac{b}{a}\cdot \frac{a}{b}C_b=C_a+C_b.$$
In particular, if $C_a$ or $C_b$ is attained, we have that $C_{a+b}<\frac{a+b}{a}C_a$ or
$C_a<\frac{a}{b}C_b$, and thus
$C_{a+b}<C_a+C_b$. We finish the proof of (2).
\ep

\bl\lab{lemma:Lions}
Let $\{u_n\}\subset S_a$ be a bounded minimizing sequence of $J[u]$ respect to $C_a<0$. Then
$$\delta:=\liminf_{n\rightarrow \infty}\sup_{z\in \R^N}\int_{B(z,1)}|u_n|^2dx>0.$$
\el
\bp
We argue by contradiction and assume that $\delta=0$, then by \cite[Lemma I.1]{Lions1984}, we have that, up to a subsequence,
$$|u_n|_p^p\rightarrow 0\;\hbox{as}\;n\rightarrow \infty, \forall p\in (2,2^*).$$
Thus $$\int_{\R^N}G(u_n)\rightarrow 0\;\hbox{as}\;n\rightarrow \infty.$$
Since $\delta=0$, one can see that
$$\int_{\R^N}V(x)u_n^2dx=o(1),$$
So
\begin{align*}
0>C_a=&J[u_n]+o(1)=\frac{1}{2}|\nabla u_n|_2^2+o(1)\geq 0,
\end{align*}
a contradiction.
\ep

\bl\lab{lemma:Lions2}
Let $\{u_n\}\subset S_a$ be a bounded minimizing sequence of $J[u]$ respect to $C_a$. If $C_a<E_a$, then there exists some $0\neq u\in H^1(\R^N)$ such that
$u_n\rightharpoonup u$ in $H^1(\R^N)$ along a subsequence.
\el
\bp
Since $u_n$ is bounded, up to a subsequence, $u_n\rightharpoonup u$ for some $u\in H^1(\R^N)$. Let $\{y_n\}\subset \R^N$ be a sequence such that
$$\int_{B(y_n, 1)}u_n^2 dx=\sup_{y\in \R^N}\int_{B(y,1)}u_n^2dx,$$
and put $\tilde{u}_n(x):=u_n(\cdot-y_n)$. Under the assumption $C_a<0$, by Lemma \ref{lemma:Lions}, we have that
$$\liminf_{n\rightarrow \infty}\int_{B(0,1)}\tilde{u}_n^2dx\geq \delta>0.$$
So if $u=0$, then we have that $\{y_n\}$ is unbounded.
Noting that $\tilde{u}_n\in S_a$,
\begin{align*}
C_a=&J[u_n]+o(1)=I[u_n]+\frac{1}{2}\int_{\R^N}V(x)u_n^2dx\\
=&I[\tilde{u}_n]+\frac{1}{2}\int_{\R^N} V(x-y_n)\tilde{u}_n^2 dx\\
=&I[\tilde{u}_n]+o(1)\geq E_a,
\end{align*}
a contradiction.
\ep

\bc\lab{cro:compact}
Under the assumptions of Lemma \ref{lemma:decreasing},
any minimizing sequence $\{u_n\}\subset S_a$ with $J[u_n]\rightarrow C_a$ and $C_a<E_a$, possesses a convergent subsequence in $H^1(\R^N)$.
\ec
\bp
By Lemma \ref{lemma:bounded-bounded}, $\{u_n\}$ is bounded in $H^1(\R^N)$. Then by Lemma \ref{lemma:Lions2}, $u_n\rightharpoonup u\neq 0$ in $H^1(\R^N)$, up to a subsequence.
We claim that $u_n\rightarrow u$ in $H^1(\R^N)$. If not, we set $v_n:=u_n-u$, then $v_n\rightharpoonup 0$ but $v_n\not\rightarrow 0$ in $H^1(\R^N)$.
Put $|u|_2^2:=\sigma$, then we see that $0<\sigma<a$ and $|v_n|_2^2=a-\sigma+o(1)$.
Then we have
\begin{align*}
C_a=&J[u_n]+o(1)=J[u]+I[v_n]+o(1).
\end{align*}
If $u$ is not a global minimizer of $J[u]$ respect to $C_\sigma$, then we have that
$$C_a>C_\sigma+I[v_n]+o(1)\geq C_\sigma+E_{a-\sigma}\geq C_\sigma+C_{a-\sigma}\geq C_a,$$
a contradiction.
If $u$ is a global minimizer of $J[u]$ respect to $C_\sigma$, then we also have
$$C_a=C_\sigma+I[v_n]+o(1)\geq C_\sigma+E_{a-\sigma}\geq C_\sigma+C_{a-\sigma}.$$
However, in this case, since $J[u]$ has global minimizer $u$ respect to $C_\sigma$, we have
$$C_{\sigma}+C_{a-\sigma}>C_a.$$
We also obtain a contradiction. Hence, we obtain that $\sigma=a$ and $J[u]=C_a$. Furthermore, we have that $u_n\rightarrow u$ in $L^2(\R^N)$. And then by Gagliardo-Nirenberg inequality, we obtain that $u_n\rightarrow u$ in $L^q(\R^N)$ for any $q\in (2,2^*)$, and then $\displaystyle \int_{\R^N}G(u_n)dx\rightarrow \int_{\R^N}G(u)dx$ as $n\rightarrow \infty$. Using the identity in Lemma \ref{lemma:Pohozaev}, we finally obtain $u_n\rightarrow u$ in $H^1(\R^N)$.
\ep

\bl\lab{lemma:not-attained}
There exists some $a_0\geq 0$ such that $C_a<0$ for $a>a_0$ and it is attained. While $C_a\equiv 0$ for $0<a<a_0$ and it is not attained.
\el
\bp
If $E_a=0$, then $C_a<0$ implies that $C_a<E_a$, then by Lemma \ref{cro:compact}, $C_a$ is attained. If $E_a<0$, then by \cite[Theorem 1.1]{Shibata2014}, $E_a$ is attained. Further, by Lemma \ref{lemma:decreasing}-(i), we see $C_a<E_a$. Then by Lemma \ref{cro:compact} again, $C_a$ is also attained.
\vskip0.11in
If the sharp number $a_0=0$, we are done. If the sharp number $a_0>0$, i.e., $C_{a_0}=0$ then by Lemma \ref{lemma:decreasing}-(1)and (3), we see that $C_a\equiv 0$ for $0<a<a_0$ and in view of Lemma \ref{lemma:decreasing}-(2) $C_a$ is not attained.
\ep

\s{Proofs  of Theorems~\ref{Main-th1}, \ref{Main-th2} and   \ref{Main-th3}}\lab{sec:proofs}
\begin{altproof}{Theorem~\ref{Main-th1}}
By Lemma \ref{lemma:J-bounded-blow}, we see that the functional $J$ is bounded from below constrained on $S_a$, i.e., $C_a>-\infty$. Let $\{u_n\}\subset S_a$ be
a minimizing sequence, that is,
$$J[u_n]\rightarrow C_a.$$
By Theorem \ref{thm:compact-infinity}, there exists some $u_0\in \mathcal{H}$ such that
$u_n\rightarrow u_0$ in both $\mathcal{H}$ and $H^1(\R^N)$. Furthermore, $u_0$ is a minimizer, i.e., $u_0\in S_a$ and $J[u_0]=C_a$.
Then there exists some $\lambda\in \R$ such that
$$-\Delta u_0+V(x)u_0-g(u_0)=-\lambda u_0.$$
\end{altproof}

\begin{altproof}{Theorem~\ref{Main-th2}}
If $E_a=0$, then $C_a<0$ implies that $C_a<E_a$, then by Lemma \ref{cro:compact}, $C_a$ is attained. If $E_a<0$, then by \cite[Theorem 1.1]{Shibata2014}, $E_a$ is attained. And by Lemma \ref{lemma:decreasing}-(i), we see $C_a<E_a$. Then by Lemma \ref{cro:compact} again, $C_a$ is also attained.

If the sharp number $a_0=0$, we are done. If the sharp number $a_0>0$, i.e., $C_{a_0}=0$. Then the desired conclusion follows from Lemma \ref{lemma:not-attained}.
 In this case, we have that
$$-\Delta u+V(x)u-g(u)=-\lambda u.$$
Noting that $C_a<0$, by $(G5)$ we have that
$$\lambda a>\lambda a+2C_a=\int_{\R^N}[g(u)u-2G(u)]dx\geq 0.$$
Hence, we have that $\lambda>0$.
\end{altproof}

\begin{altproof}{Theorem~\ref{Main-th3}}
Under the assumption (i), by Remark \ref{remark:zy2} and \cite[Theorem 1.3]{Shibata2014}, we obtain that $a_0\leq \alpha_0=0$.

Under the assumption $(ii)$, by \cite[Theorem B]{Ikoma-Miyamoto2020}, we also have that $a_0=0$.
\end{altproof}

\noindent
{\bf Acknowledgements}\\
The authors thank Prof. Louis Jeanjean very much for the valuable comments and suggestions when preparing the paper.

\vskip0.26in


\end{document}